\pdfoutput=1
\documentclass{scrartcl}

\usepackage{amssymb,amsmath,bm,graphicx,graphbox,
  dsfont,xcolor,float,booktabs,hyperref,units,mathtools,subcaption,tikz}
\usepackage[labelfont=bf]{caption}
\usetikzlibrary{decorations.pathreplacing,fadings}
\usetikzlibrary{shapes.geometric,positioning}
\usepackage[Algorithm]{algorithm}
\usepackage{algpseudocode}
\usepackage{setspace}
\usepackage{mycommands}

\newcommand{\blue}[1]{{\color{blue} \noindent #1}}

\title{Structure Preservation for the Deep Neural Network Multigrid
  Solver}
\author{Nils Margenberg\footnotemark[2]
        \and Christian Lessig\footnotemark[3],
        \and Thomas Richter\footnotemark[4]}

\begin{document}

\maketitle

\begin{abstract}
  The simulation of partial differential equations is a central
  subject of numerical analysis and an indispensable tool in
  science, engineering and related fields. Existing approaches, such
  as finite elements, provide (highly) efficient tools but deep neural
  network-based techniques emerged in the last few years as an
  alternative with very promising results. We investigate the
  combination of both approaches for the approximation of the Navier-Stokes
  equations and to what extent structural properties such
  as divergence freedom can and should be respected.
  Our work is based on DNN-MG, a deep neural
  network multigrid technique, that we introduced recently
  and which uses a neural network to represent fine grid
  fluctuations not resolved by
  a geometric multigrid finite element solver.
  Although DNN-MG provides solutions with very good
  accuracy and is computationally highly efficient,  we noticed that the
  neural network-based corrections substantially violate the
  divergence freedom of the velocity vector field.
  In this contribution, we discuss these findings and analyze three
  approaches to address the problem: a penalty term to encourage
  divergence freedom of the network output; a penalty term for the corrected velocity field;
  and a network that learns the stream function,
  i.e.\ the scalar potential of the divergence free velocity vector field and which
  hence yields by construction divergence free corrections.
  Our experimental results show that the third approach based on
  the stream function outperforms the other two and not only improves the
  divergence freedom but in particular also the overall fidelity of the simulation.
\end{abstract}

\renewcommand{\thefootnote}{\fnsymbol{footnote}}

\footnotetext[2]{Helmut Schmidt University,
Holstenhofweg 85, 22043 Hamburg,
Germany,
\texttt{margenbn@hsu-hh.de}}
\footnotetext[3]{University of Magdeburg,
Institute for Simulation and Graphics,
Universit\"atsplatz 2,
39104 Magdeburg,
Germany,
\texttt{christian.lessig@ovgu.de}}
\footnotetext[4]{University of Magdeburg,
Institute for Analysis and Numerics,
Universit\"atsplatz 2,
39104 Magdeburg,
Germany,
\texttt{thomas.richter@ovgu.de}}

\section{Introduction}
\label{sec:intro}

The tremendous progress on deep learning in the last decade for applications such as image recognition and machine translation~\cite{LeCun2015} has motivated their use also for the simulation of physical systems.
Despite promising results that have already been obtained on this, many fundamental questions remain open.
One of them is how physical constraints, such as energy and momentum conservation or divergence freedom, can be integrated into a neural network-based simulation---and if an explicit consideration of this is necessary or meaningful.

In the literature, two principle approaches to enforce physical constraints have been proposed.
In the first one, they are imposed by construction, e.g.\ through an appropriate network design or by representing the physical quantities so that their properties and relationship have to be satisfied.
Examples are symplectic neural network~\cite{Mattheakis2019,Chen2019} that impose symplecticity by construction or the representation of fluid velocity fields in a divergence free basis with learning being performed on the basis function coefficients~\cite{Sapsis2019,Wan2020}.
The second and more common approach to physical constraints in neural networks is to encourage that these hold but not strictly enforce them.
Such a ``soft'' formulation is typically accomplished using a penalty term in the loss function but can also be realized, e.g., with a Bayesian formulation with the the physical constraints (including an entire partial differential equation) as priors.
The approach also naturally accommodates real-world data where properties such as energy conservation almost never perfectly hold (e.g.\ because a system is not closed).
However, to what extent a constraint is satisfied after training is difficult to control through a penalty term and its addition to the loss function increases the training complexity and can hinder convergence.
%With constraints enforced by construction the network is restricted to a valid parameter domain which typically simplifies training.

In this work, we explore how divergence freedom of the fluid velocity vector field can be ensured when a neural network is used in the simulation of the incompressible Navier-Stokes equations.
In particular, we build on the recently introduced Deep Neural Network Multigrid Solver (DNN-MG)~\cite{Margenberg2020} that tightly integrates a neural network into a multigrid finite element solver to replace the computations on one or multiple finest mesh levels by a neural network-based prediction.
In~\cite{Margenberg2020}, the efficiency of DNN-MG for the simulation of the incompressible, instationary Navier-Stokes equations was demonstrated.
However, the predicted velocity vector field is not divergence free.
We hence study in this paper three different approaches to encourage or enforce divergence freedom within the DNN-multigrid formulation.
First, we modify the network's loss function and include a penalty term to limit the divergence of the predicted velocity field.
Second, we penalize the divergence of the corrected velocity field, which has the potential to also reduce the divergence residual of the finite element velocity field.
And third, we design a neural network that learns the stream function, i.e.\ the potential associated with the velocity vector field, and hence provides by construction divergence free corrections.

Our results demonstrate that encouraging or enforcing divergence freedom improves the overall quality of the simulations, measured, e.g., by the drag and lift functionals, although the three approaches we developed lead only to rather modest reductions in the divergence.
For the two approaches based on a penalty term, this results from a delicate balance between prediction accuracy and divergence freedom that needs to be attained but that it is difficult to satisfy in practice.
The approach based on the stream function provides overall the best results and achieves lift and drag functionals close to a fine reference solution.
At the same time, the corresponding network is more efficient to train.
We attribute this improved performance to a formulation that better respects the intrinsic structure of the partial differential equations.

The remainder of the paper is structured as follows.
After discussing related work, we recall the DNN-MG solver in Sec.~\ref{sec:dnnmg}.
In Sec.~\ref{sec:structure} we describe our approaches to ensure divergence freedom within DNN-MG and present and discuss our numerical results.
We conclude in Sec.~\ref{sec:conclusion} where we also provide an outlook for possible directions for future work.

\section{Related Work}
\label{sec:related}

%The advent of deep neural networks is yielding tremendous progress on a wide range of problems in machine translation, computer vision, and related fields~\cite{LeCun2015}. 
Recently, there has been an increasing interest to employ (deep) neural networks for the simulation of physical systems.
In the following, we will discuss existing approaches with an emphasis on those that consider partial differential equations and that are most pertinent to our own work.

%\paragraph{Neural network architectures for partial differential equations}
Neural network-based techniques for partial differential equations can, broadly speaking, be classified into three categories. 
The first one, e.g.\ the DeepRitz method by E and Yu~\cite{E2018} or the physics-informed neural networks by Raissi, Karniadakis and co-workers~\cite{Raissi2019b,Lu2020b}, use neural networks in place of a basis or frame to represent solutions.
Training then amounts to fitting the network representation to the initial or boundary value problem.
This has hence to be repeated for each problem instance but the approach yields typically a mesh free solution that scales well to high dimensions~\cite{E2018} and the methodology applies to a wide range of partial differential equations. 
The second approach takes inspiration from classical time stepping schemes and considers the solution of the PDE as the problem of predicting the next state from the previous ones. 
The well developed network architectures for sequential data can hence be employed, e.g.\ recurrent neural networks (RNN) with Long-Short-Term-Memory units (LSTM)~\cite{Hochreiter1997} or Gated Recurrent Units (GRUs)~\cite{Cho2014}, or temporal convolutional networks, cf.~\cite{Bai2018}.
Examples of this approach are~\cite{Hu2020,Stevens2020,Wiewel2018} and the work by Han, Jentzen, and E~\cite{Han2018} who use a formulation of semi-linear, parabolic PDEs as backward stochastic differential equations and then realize each time step of these as a layer in a deep neural network.
Since we combine a classical finite element based simulation, which relies on a time stepping scheme, with a neural network, our work also falls into this category.
The third class of neural network-based techniques for partial differential equations uses these to represent solution operators.
They hence use the initial condition as well as parameters on which the solution depends as input to the neural network and it then maps these to the solution.
For example, in~\cite{Bhattacharya2020} the solution operator for parametric partial differential equations is learned based on a discretization obtained via model reduction and~\cite{Li2020b,Li2020} learn mappings between infinite dimensional spaces that can be solution operators.
Lu, Jin, Karniadakis~\cite{Lu2020} learn nonlinear operators from data, based on the result that any such operator can be approximated by a neural network~\cite{Chen1995}.

%The most common architecture are recurrent neural networks (RNN), e.g.\ with Long-Short-Term-Memory units (LSTM)~\cite{Hochreiter1997} to avoid the vanishing gradient problem. 
%Because of the prevalence of sequential data, LSTMs led to many variations over the years.
%One are Gated Recurrent Units (GRUs) which, while not as powerful as LSTMs, provide often similar performance and are easier to train~\cite{Cho2014}.
%Our work fo

%We will use these in our work.

%An alternative to LSTM-type architectures are temporal convolutional neural networks (TCNs) where the temporal dependence is modeled using a causal (i.e.\ one-sided) convolution in the time variable.
%Bai et al.~\cite{Bai2018} demonstrated that this can provide substantial improvements in prediction accuracy over recurrent neural networks with memory.
%A different approach to modeling temporal dependencies was developed by Voelker et al.~\cite{Voelker2019} who represented the time domain in Legendre polynomials and learned basis function coefficients.
%Through this, they obtained a neural network model with continuous time.
%A related idea are the neural ordinary differential equations architectures by Chen et al.~\cite{Chen2018} where the depth of the network is continuous.

%\paragraph{Deep learning of partial differential equations}
%Approaches for a direct description of the time evolution of partial differential equations using neural networks go back more than 20 years, e.g.~\cite{Lagaris1998}.
%New impetus for this direction arose recently through the advances made with deep neural networks for other problems.
For flow problems, Raissi et al.~\cite{Raissi2020} demonstrated that the velocity and pressure of the Navier-Stokes equations can be learned to good approximation only from observing passively advected tracer particles. 
Nabian and Meidani~\cite{Nabian2018}, Yang and Perdikaris~\cite{Yang2018}, and Raissi and Karniadakis~\cite{Raissi2018} exploited that not only observations are available but also a known analytic model.
These authors hence also use the deviation of the network prediction from the model as an additional penalty term in the training loss.
Kasim et al.~\cite{Kasim2020} recently demonstrated neural network-based simulations for a broad range of applications, including fluid dynamics, by also optimizing the network architecture itself during training.
Eichinger, Heinlein and Klawonn~\cite{EichingerHeinleinKlawonn2020} use convolutional neural networks and techniques from image processing to learn flow patterns for the Navier-Stokes flow around objects of different shape. 
%The present project is complementary since we would like to use the neural networks in particular also to represent phenomena for which no model exists.

%\paragraph{Deep learning assisted of partial differential equations}
Next to the above direct approaches to the simulation of partial differential equations using neural networks, there have been different attempts to integrate these into existing numerical formulations.
As in our work, the objective is then to combine the benefits of neural networks and classical methods and obtain techniques that would be difficult or impossible with either approach alone.
For elasticity, neural network-based representations for constitutive relations were learned in ~\cite{Tartakovsky2018,Berg2019,Rudy2019} and then used in classical simulations, e.g.\ based on finite elements.
%i.e.\ type B problem in the jargon of~\cite{Abdulle2012}. 
Wiewel et al.~\cite{Wiewel2018} presented a simulation of the Navier-Stokes equation where an LSTM-based deep neural network is used to predict the pressure correction within a classical finite element simulation. 
In~\cite{Long2019} a convolutional neural network is interpreted as a finite difference discretization of a PDE and this was also proposed in~\cite{Ruthotto2020}.
Bar-Sinai, Hoyer, Hickeney and Brenner~\cite{BarSinai2019} learn general sub-grid scale models, used, e.g., with finite difference or finite volume discretization.
In the context of large eddy-type simulations, the learning of closure models has also been considered, e.g.\ in~\cite{Ling2016}.
Wan, Dodov, Lessig, Dijkstra, and Sapsis~\cite{Wan2020} combine a wavelet-based discretization and perform learning in the coefficient space where space-frequency correlations are resolved but, so far, without simulation.
Recently, Stevens and Colonius~\cite{Stevens2020} enriched finite difference and finite volume simulations with a deep neural network to resolve sharp discontinuities (as in Burger's equations) and to deal with PDEs with chaotic behavior.

%An approach based on the splitting method for high dimensional PDEs was developed by Beck et al.~\cite{Beck2019} where one obtains a set of smaller learning problems that are easier to treat than a single large one.
%Wan and Sapsis~\cite{Wan2018} describe the motion of finite size particles in fluid flows by combining analytic arguments with a neural network that models the discrepancy between the (highly) idealized analytic model and real world observations.

To ensure that neural network-based predictions respect physical invariants, such as energy conservation or divergence freedom, network architectures tailored towards physical simulations have been proposed.
Next to an implicit enforcement as in Physics Informed Neural Network (PINN)~\cite{Lu2020b}, various authors proposed architecture where desired properties hold by construction.
Finizi, Wang and Wilson~\cite{Finzi2020} and Cranmer et al.~\cite{Cranmer2020}, e.g., learn the dynamics of classical mechanical systems by directly learning the Hamiltonian or Lagrangian functions.
In related work,~\cite{Mattheakis2019,Chen2019,Jin2020} developed neural networks that incorporate the symplectic structure of Hamiltonian mechanics and they demonstrate that this improves generalization and prediction accuracy.
Based on the work by Bruna and Mallat~\cite{Bruna2013}, Bietti and Mairal~\cite{Bietti2019} proposed neural network architectures that are invariant under group transformation, including diffeomorphisms.
A neural network that is directly defined on a Lie group, and hence also enjoys group equivariance, was proposed by Bekker~\cite{Bekkers2020}.

\section{Deep Neural Network Multigrid Solver}
\label{sec:dnnmg}

In this section we summarize the Deep Neural Network Multigrid Solver (DNN-MG) introduced in~\cite{Margenberg2020}.
For a detailed description the reader is referred to the original paper.

%%%%%%%%%%%%%%%%%%%%%%%%%%%%%%%%%%%%%%%%%%%%%%%%%%%%%%%%%%%%%%%%%%%%%%%%%%%%%%%%
\subsection{Finite Element Discretization of the Incompressible Navier-Stokes equations}

We want to solve the incompressible, instationary Navier-Stokes equations given by
\begin{equation}
  \begin{alignedat}{2}
  \label{eq:nsstrong}
    \partial_t v + (v\cdot \nabla)v - \frac{1}{\mathrm{Re}}\Delta v
    +\nabla p &= f \quad &&\text{on } [0,\,T] \times \Omega\\
    \nabla \cdot v &= 0 \quad &&\text{on } [0,\,T] \times \Omega,
  \end{alignedat}
\end{equation}
where $v\colon [0,\,T]\times \Omega \to \R^2$ is the velocity, $p\colon [0,\,T]\times \Omega \to \R$ the pressure, $\mathrm{Re}>0$ the Reynolds number, and $f$ an external force.
The initial and boundary conditions are given by
\begin{equation}
  \label{eq:boundary}
  \begin{alignedat}{2}
    v(0,\,\cdot) &= v_0(\cdot)\quad &&\text{on }\Omega\\
    v &= v^D \quad &&\text{on } [0,\,T] \times \Gamma^D\\
    \frac{1}{\mathrm{Re}}(\vec n\cdot\nabla)v - p\vec n &=0 \quad &&\text{in } [0,\,T] \times \Gamma^N ,
  \end{alignedat}
\end{equation}
where $\vec n$ denotes the outward facing unit normal on the boundary $\partial\Omega$ of the domain. On the outflow boundary $\Gamma^N$ we consider the do-nothing outflow condition~\cite{HeywoodRannacherTurek1992}, which is well established to model artificial boundaries; see the discussion in~\cite{BraackMucha2014} on possible shortcomings and a variation that is able to guarantee uniqueness of a solution.

We discretize Eq.~\ref{eq:nsstrong} using a weak finite element formulation with $\smash{v_h,\,\phi_h \in V_h = [W_h^{(2)}]^d}$ and $\smash{p_h,\,\xi_h \in L_h = W_h^{(2)}}$ where $\phi_h$ and $\xi_h$ are test functions and $\smash{W_h^{(r)}}$ is the space of continuous functions which are polynomials of degree $r$ on each mesh element $T \in\Omega_h$ in the mesh domain $\Omega_h$.
The resulting equal order finite element pair $V_h\times L_h$ does not fulfill the inf-sup condition.
We hence use stabilization terms of local projection type~\cite{BeckerBraack2001} with parameter $\alpha_T = \alpha_0\cdot  \mathrm{Re} \cdot h_T^2$ and projection $\smash{\pi_h:W_h^{(2)}\to W_h^{(1)}}$ into the space of linear polynomials.

With the second order Crank-Nicolson method for time discretization, the solution of Eq.~\ref{eq:nsstrong} subject to Eq.~\ref{eq:boundary} at each time step amounts to determining the state $x_n=(v_n,p_n)$ such that
\begin{equation}
  \label{eq:nonlinearshort}
    {\cal A}_n(x_n) = f_n
  \end{equation}
where
  \begin{equation}\label{varform}
  \begin{aligned}
    [{\cal A}_n(x)]_i &\coloneqq
    (\nabla \cdot v_{n}, \, \xi_h^i)+
    \sum_{T\in\Omega_h}\alpha_T (\nabla (p_n-\pi_h p_n),\nabla (\xi_h^i-\pi_h \xi_h^i))\\
    & \qquad +\frac{1}{k}(v_n,\,\phi_h^i)\,
    +{}\frac{1}{2} (v_n\cdot \nabla v_n,\,\phi_h^i)
    +\frac{1}{2\mathrm{Re}}(\nabla v_n,\,\nabla \phi_h^i)
    -(p_n,\,\nabla \cdot \phi_h^i)
    \\[6pt]
    [f_n]_i &\coloneqq  \frac{1}{k}(v_{n-1},\, \phi_h^i)
    +\frac{1}{2}(f_n,\phi_h^i)
    +\frac{1}{2}(f_{n-1},\,\phi_h^i)
    \\
    &\qquad
    -\frac{1}{2}(v_{n-1}\cdot\nabla v_{n-1},\,\phi_h^i)
    - \frac{1}{2\mathrm{Re}}(\nabla v_{n-1},\,\nabla \phi_h^i),
  \end{aligned}
\end{equation}
for all test functions $\phi_h^i$ and $\xi_h^i$.
Eq.~\ref{eq:nonlinearshort} is a large nonlinear system of algebraic equations that is solved by Newton's method based on the initial guess $x_n^{(0)}=(v_{n-1},p_{n-1})$ and the iteration
\begin{equation}
  \label{eq:newton}
  {\cal A}_n'(x_n^{(l-1)}) \, w^{(l)} = f_n-{\cal A}_n(x_n^{(l-1)}), \quad \quad x_n^{(l)}=x_n^{(l-1)}+w^{(l)}
\end{equation}
for $l=1,2,\dots$
Here we denote by ${\cal A}'(x^{(l-1)})$ the Jacobian of ${\cal A}$ at $x^{(l-1)}$, which for the problem at hand can be computed analytically, cf.~\cite[Sec. 4.4.2]{Richter2017}.

Each Newton step requires the solution of the linear system in Eq.~\ref{eq:newton}, where the system matrix ${\cal A}'(x^{(l-1)})$ is sparse but non-symmetric and not definite due to the saddle point structure of the underlying Navier-Stokes equation. To approximate the solution to Eq.~\ref{eq:newton} with optimal robustness, we employ the generalized minimal residual method (GMRES) introduced by Saad~\cite{Saad1996}. The convergence is accelerated through a preconditioner that is realized by a single sweep of a geometric multigrid solver.

%\begin{algorithm}[t]
%  \caption{The geometric multigrid method for the solution of the linear system $A_L \, x_L = b_L$ given on a finest level $L$.
%    If all high frequency errors of $S$ are smoothed
%    at a constant rate, the method achieves optimal complexity
%    $O(n)$. The  multigrid solver is initiated on the finest mesh level $L$ and then used recursively. In our work it is employed for the solution of the linear system in Eq.~\ref{eq:newton} required for every Newton step for the resolution of the nonlinear problem in Eq.~\ref{eq:nonlinearshort}.}
%  \label{alg:geomg}
%  \begin{algorithmic}[1]
%    \Procedure{multigrid}{$l,\;A_l,\;b_l,x_l$}
%    \State{$s_l \leftarrow S(A_l,\;b_l,\;x_l)$}\Comment{Smoothing}
%    \State{$r_l \leftarrow b_l - A_l s_l$}\Comment{Residual}
%    \State{$r_{l-1} \leftarrow \RS(l,\,r_l)$}\Comment{Restriction}
%    \If{$l-1=0$}\Comment{Coarse-grid solution}
%    \State{$c_0 \leftarrow A_0^{-1}r_0$}\Comment{Direct solution}
%    \Else{}
%    \State{$c_{l-1} \leftarrow \Call{multigrid}{l-1,\,A_{l-1},\,r_{l-1},\,0}$}
%    \EndIf{}
%    \State{$x'_l \leftarrow s_l + \PL(l,\,c_{l-1})$}\Comment{Prolongation}
%    \State{$s'_l \leftarrow S(A_l,\;b_l,\;x'_l)$}\Comment{Smoothing}
%    \State{\Return{$s'_l$}}
%    \EndProcedure
%  \end{algorithmic}
%\end{algorithm}

The geometric multigrid method is based on a hierarchical approximation of a linear system on a sequence of finite element spaces
$V_0\subset V_1\subset \cdots \subset V_L$ defined over a
hierarchy of meshes $\Omega_0,\dots,\,\Omega_L=\Omega_h$.
Instead of treating the linear system on its domain, given by the finest
mesh level $L$ (as one would do in a traditional solver)
the high frequencies there are smoothed
and the remaining errors are treated on lower levels.
%Usually, smoothing is
%applied at the beginning %(Algo.~\ref{alg:geomg}, line 2)
%and at the end %(Algo.~\ref{alg:geomg}, line 11)
%of each iteration and in between
%the residual is computed %(Algo.~\ref{alg:geomg}, line 3)
%and restricted to the next
%coarse level $L-1$. %(Algo.~\ref{alg:geomg}, line 4).
This is accomplished by computing the residual
on level $L$ %(Algo.~\ref{alg:geomg}, line 3)
and restricting it to the next
coarse level $L-1$. %(Algo.~\ref{alg:geomg}, line 4).
There, the process is repeated and
this is performed recursively %(Algo.~\ref{alg:geomg}, line 8)
until a coarsest level has been reached where
a direct solver is employed. %(Algo.~\ref{alg:geomg}, line 6).
The updates from the coarser levels
are then prolongated back to level $L$ and the coarse-to-fine and
fine-to-coarse iteration is repeated until a prescribed
error tolerance has been met. %(Algo.~\ref{alg:geomg}, line 10).
The mesh transfer from fine to coarse is accomplished with
$L^2$-projections, known as restrictions, and from coarse to fine
with interpolations, known as prolongations.
The smoothing that preceded the restrictions is realized with a
simple smoothing operator $S(A_l,b_l,x_l)$ that yields an approximate
solution of the linear system $A_l \, x_l = b_l$, i.e.
$S(A_l,b_l,x_k) \approx A_l^{-1}b_l$
and that aims to quickly reduce all high frequency components of the residual $b_l-A_lx_l$.
%If this happens on all levels $l$ at a constant rate, then the geometric multigrid method achieves the optimal complexity $\mathcal{O}(n)$ for the solution of the linear system.
In our implementation we use a simple iteration of Vanka type~\cite{Vanka1985}, which allows for easy parallelization and gives very good performance with less than 5 pre- and post-smoothing steps~\cite{FailerRichter2020}.
%The idea of the Vanka type smoother is to exactly solve small subproblems and to replace the inverse of $A_l^{-1}$ by
%\begin{equation}\label{vanka0}
%  {\cal V}_l(A_l) =\sum_{T\in\Omega_l}R_T^T [R_T A_lR_T^T]^{-1}R_T
%\end{equation}
%where $R_T$ is the restriction to those nodes that belong to one mesh element $T\in\Omega_l$ and $R_T^T$ its transpose. Considering piecewise quadratic finite elements, 9 nodes are attached to each element such that the local matrices to be inverted have the dimension $R_TA_lR_T^T \in\mathds{R}^{27\times 27}$, since each node comprises the scalar pressure and two velocity components.
%  With a damping parameter $\omega\in (0,1]$, one iteration of the Vanka smoother is then given by
%  \begin{equation}\label{vanka}
%    S(A_l,b_l,x_k)=x_k + \omega {\cal V}_l(A_l) (b_l-A_lx_k).
%  \end{equation}

%%%%%%%%%%%%%%%%%%%%%%%%%%%%%%%%%%%%%%%%%%%%%%%%%%%%%%%%%%%%%%%%%%%%%%%%%%%%%%%%
\subsection{Deep Neural Network Multigrid Solver}

\begin{algorithm}[t]
  \setstretch{1.2}
  \caption{DNN-MG for the solution of the Navier-Stokes equations. Lines 6-9 (blue) provide the modifcations of the DNN-MG method compared to a classical Newton-Krylow simulation with geometric multigrid preconditioning.}
  \label{alg:dnnmg}
  \begin{algorithmic}[1]
    \For{all time steps $n$}
%    \State{$b_{\scriptscriptstyle L+1}^n \leftarrow \textbf{Rhs}(v^{n-1}_{\scriptscriptstyle L+1})$}
%    \Comment{Compute rhs of Eq.~\ref{eq:4cranknicholson3}}
%    \State{$b_{\scriptscriptstyle L}^n \leftarrow \mathcal{R}(b_{\scriptscriptstyle L+1}^n)$}
%    \Comment{Restrict rhs to $L$}
    \While{not converged}\Comment{Newton's method in Eq.~\ref{eq:newton}}
      \State{$\delta x_i \leftarrow$ \Call{multigrid}{$L,\,A_{L}^n,\,b_{L}^n,\,\delta x_i$}}\Comment{Geometric multigrid with $x = (p_n^L,\,v_n^L)$}
      \State{$x_{i+1} \leftarrow x_i + \epsilon \, \delta x_i$}
    \EndWhile{}
    \State{\blue{$\tilde{v}_n^{\scriptscriptstyle L+1} \leftarrow \mathcal{P}(v_n^{\scriptscriptstyle L}) $}}\Comment{\blue{Prolongation on level $L+1$}}
    \State{\blue{$d_n^{\scriptscriptstyle L+1} \leftarrow \mathcal{N}(\tilde{v}_n^{\scriptscriptstyle L+1},\,\Omega_L,\,\Omega_{L+1})$}}\Comment{\blue{Prediction of velocity correction}}
    \State{\blue{$b_{n+1}^{\scriptscriptstyle L+1} \leftarrow \mathrm{Rhs}(\tilde{v}_n^{\scriptscriptstyle L+1} + d_n^{\scriptscriptstyle L+1},f_n,f_{n+1})$}}\Comment{\blue{Set up rhs of Eq.~\ref{eq:nonlinearshort} for next time step}}
    \State{\blue{$b_{n+1}^{\scriptscriptstyle L} \leftarrow \mathcal{R}(b_{n+1}^{\scriptscriptstyle L+1})$}}\Comment{\blue{Restriction of rhs to level $L$}}
    \EndFor{}
  \end{algorithmic}
\end{algorithm}

The solution of the Newton iteration in Eq.~\ref{eq:newton} using GMRES and the geometric multi-grid method is highly efficient.
However, since one has multiple GMRES steps for every step of the Newton method and for each one also one up- and down-sweep of the multigrid method, there is still a substantial amount of computations.
The majority is thereby required for the finest mesh level (typically about three quarters in the multigrid method) so that a refinement of the mesh, or, equivalently, the use of an additional mesh level $L+1$, incurs a substantial increase in the computational costs.
The Deep Neural Network Multigrid Solver (DNN-MG) addresses this bottleneck and replaces computations on one or multiple finest mesh levels with a neural network-based prediction.
An overview of DNN-MG is provided in Algo.~\ref{alg:dnnmg}.

%  \begin{figure}[t]
%    \begin{center}
%      \includegraphics[width=0.9\textwidth]{figures/patches.pdf}
%    \end{center}
%    \caption{Setup of the DNN-MG solver. The coarse mesh solution on the
%      mesh $\Omega_L$ is prolongated to the fine prediction level
%      $\Omega_{L+1}$. Here the neural network is acting locally on a
%      patch structure. We show two different patch sizes $p=0$ and $p=1$
%      that determine the input size of the neural network. In blue, we
%      indicate the degrees of freedom associated with a piecewise
%      bi-quadratic finite element space. \red{Not clear what a patch is and with what color it is indicated. Red? Then probably best to show just one patch and say so in the text. I also don't think the $p$ is needed (we didn't have it in the JCP submission) and it does not appear later.} }
%    \label{fig:patch}
%  \end{figure}

To ensure that the DNN-MG solver is efficient and flexible, its neural network is build around a few critical design principles:
\begin{enumerate}
  \item The neural network operates patch-wise, i.e.\ locally on small neighborhoods of the mesh domain $\Omega_{L+k}$. This ensures the generalizability of DNN-MG to different flow regimes and meshes.
%      A patch is a local and structured set of elements, ranging from a single element to $(2^p)^d$ adjacent elements, where $d$ is the spatial dimension and $p$ the patch-depth of the prediction layer. Compare Fig.~\ref{fig:patch}. % while enabling a highly compact network architecture.
  \item The neural network has memory. This ensures that complex flow behavior can be predicted and that corrections are coherent in time.
  \item The neural network uses the residual of the nonlinear problem in Eq.~\ref{eq:nonlinearshort} on level $L+k$ as input. This provides rich information about the required correction. % and through this again facilitates a compact network architecture.
\end{enumerate}
Together, these design principles lead to a surprisingly compact neural network architecture that at the same time generalizes well to different flow regimes and domains.
In particular, the implementation in~\cite{Margenberg2020} has just 8634 trainable parameters with an architecture based on Gated Recurrent Units.
It can hence be trained in a few hours with a small number of example flows as training set.
%As loss function for the training a simple $l_2$ loss between the predicted velocity and a precomputed reference was used.

As shown in Algo.~\ref{alg:dnnmg}, the neural network-based correction of DNN-MG is applied at the end of every time step after the Newton solver has computed an updated velocity $v_n^L$ on level $L$, cf. Eq.~\ref{eq:newton}.
To obtain the correction, $v_n^L$ is first prolongated to level $L+1$, yielding $\tilde{v}_n^{L+1} = \mathcal{P}(v_n^L)$.
The neural network is then employed individually for each patch $P_i$ (a mesh element on level $L$ in the most simple case as in~\cite{Margenberg2020}) and predicts a velocity correction $d_{n,i}^{L+1}$.
The inputs to the neural network are thereby also entirely local and include the residual of Eq.~\ref{eq:nonlinearshort} on level $L+1$, the prolongated velocity $\smash{\tilde{v}_n^{L+1}}$, geometric properties such as the patch's aspect ratio, and possibly other quantities such as the P{\'e}clet number over the patch.
With the correction $d_{n}^{L+1}$, a provisional right hand side $b_{n+1}^{\scriptscriptstyle L+1} = \mathrm{Rhs}(\tilde{v}_n^{\scriptscriptstyle L+1} + d_n^{\scriptscriptstyle L+1},f_n,f_{n+1})$ of Eq.~\ref{eq:newton} is formed on level $L+1$ and then restricted to level $L$, i.e.\ we compute $b_{n+1}^{\scriptscriptstyle L} = \mathcal{R}(b_{n+1}^{\scriptscriptstyle L+1})$. Restricting the right hand side instead of the corrected velocity $\tilde{v}_n^{L+1} + d_n^{L+1}$ is thereby essential to propagate information back to level $L$ since, ideally, $\mathcal{R}( \mathcal{P}(v_n^L) + d_n^{L+1}) = v_n^L$.
The corrected $b_{n+1}^{\scriptscriptstyle L}$ is then used in the next time step in the Newton solve, which is again improved by a neural network-based correction at the end of the time step.

%    \blue{The neural network is applied patch by patch and the global velocity update is assembled by locally averaging over all nodes are part of multiple overlapping patches, see also Fig.~\ref{fig:patch}
%    \begin{equation}\label{average}
%      {\cal N}(v^L,\Omega_L) \coloneqq
%      \sum_{P_i\in\Omega_{L+1}}A_{P_i} {\cal N}(v^L|_{P_i},P_i),
%    \end{equation}
%    where we denote by $A_{P_i}$ the averaging and prolongation operator corresponding to the patch $P_i$
%    \begin{equation}\label{average:1}
%    \Big(A_{P_i}\Big)_{kl} = \begin{cases}
%      0 & \text{global degree of freedom }k\text{ is not part of patch }P_i,\\
%      & \text{or it does not correspond to node }l\text{ in patch }P_i.\\
%      \frac{1}{e_k} & \text{global degree of freedom }k\text{ corresponds to node }l\text{ within patch }P_i\\
%      &\text{and it is part of }n_k\text{ different patches}.
%    \end{cases}
%    \end{equation}
%    The dependency of the network on the mesh $\Omega_L$, locally expressed by the geometry of the patches $P_i$, is added to allow the network to abstract information on the locally varying geometry and mesh size.}

The neural network of DNN-MG is trained using a high fidelity finite element solution obtained on the fine mesh level $L+1$ and with the loss function 
% is constructed locally on each patch
\begin{equation}\label{loss}
  \mathcal{L}(v^L,v^{L+1};d^{L+1})
  \coloneqq \sum_{n=1}^{N}\sum_{P_i\in \Omega_{L+1}}
  \big\Vert v^{L+1}_{n+1} - (v^L_n+d^{L+1}_n)\big\Vert^2_{l^2(P_i)}.
\end{equation}
Here, $N$ is the number of time-steps that make up the training data and, since DNN-MG operates strictly local over patches, for each one the loss accumulates the local residuals in the second sum.
%  $v^L_n$ is the coarse mesh finite element multigrid solution at time $t_n$ and $v^{L+1}_{n+1}$ is the high fidelity solution obtained by a resolved simulation on the fine mesh at time $t_{n+1}$. By $d^{L+1}_n={\cal N}(v^L_n)$ we denote the network output. 
Due to this strictly local structure of DNN-MG, a very large data set is generated by a single resolved simulation.

In~\cite{Margenberg2020} it was demonstrated that the DNN-MG algorithm substantially improves lift and drag functionals for the classical channel flow around a cylinder and generalized well to flow regimes not seen during training, e.g.\ a flow without an obstacle or a flow in an L-shaped domain.
%We refer to~\cite{Margenberg2020} for more details on the deep neural network multigrid solver.

\section{Enforcing Divergence Freedom in DNN-MG}
\label{sec:structure}

In our first study~\cite{Margenberg2020}, we have shown that DNN-MG is able to increase the approximation quality at very little extra cost. The error in the velocity vector field but also in derived quantities such as the drag and lift functionals of the flow around an obstacle are reduced considerably with an increase in the computation time that is orders of magnitude below those for adding an additional mesh level.
We have, however, also observed that the network is not good at ensuring a solenoidal velocity vector field. Although the loss-function aims at minimizing the distance to the high fidelity solution, which weakly satisfies divergence freeness to a certain degree, the neural network minimizes Eq.~\ref{loss} with a vector field with a substantial divergence.

\subsection{Violation of Divergence Freedom in Existing Formulation}

\begin{figure}[t]
  \includegraphics[width=0.5\textwidth]{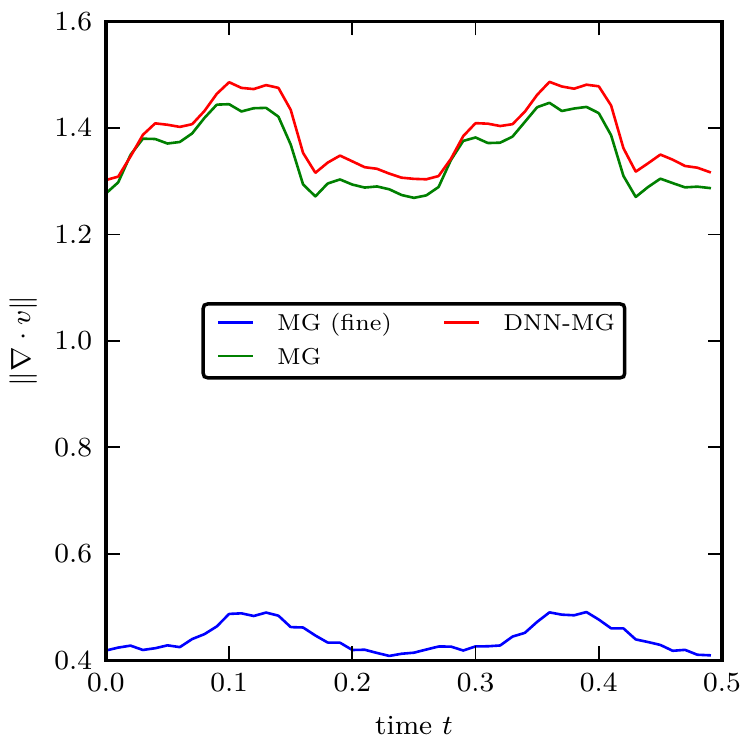}
  \includegraphics[width=0.5\textwidth]{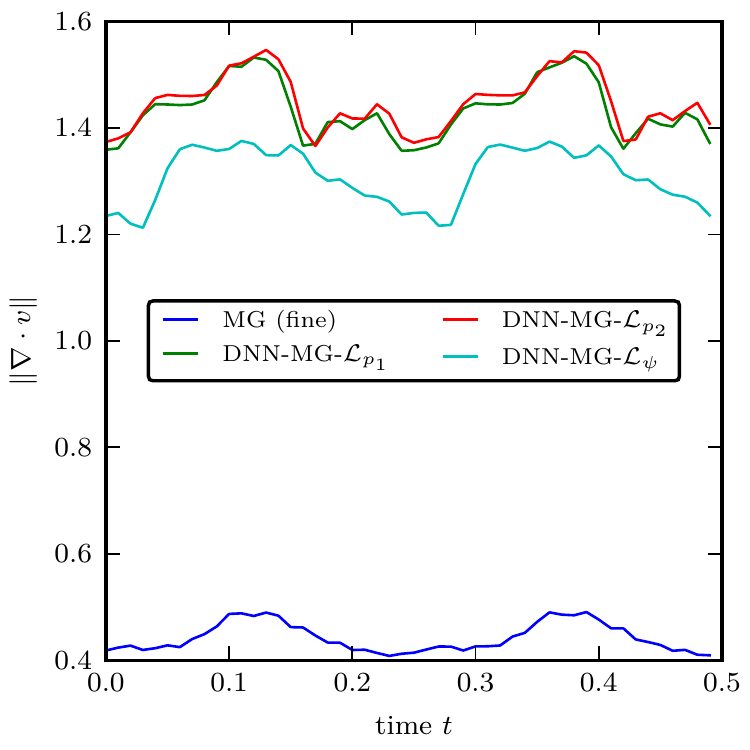}
\caption{
  Divergence $\|\nabla\cdot v\|$ for the coarse mesh (MG),
  the high fidelity (MG fine) and the deep neural network multigrid solution (DNN-MG).
  Three results of the attempts to enforce divergence freedom are shown:
  DNN-MG with the divergence of the DNN's output as penalty term (DNN-MG-$\mathcal{L}_{p_1}$),
  DNN-MG with the divergence of the corrected velocity field as penalty term
  (DNN-MG-$\mathcal{L}_{p_2}$) and enforcing divergence freedom by construction (DNN-MG-$\mathcal{L}_{\psi}$).}
  \label{fig:div}
\end{figure}

The quadratic equal-order finite element approach $(v_h,p_h)\in V_h\times L_h$ that provides the basis for the DNN-MG algorithm in our current implementation is stabilized with local projections~\cite{BeckerBraack2001}, cf.~Eq.~\ref{varform}.
The discrete divergence equation thus reads
\[
(\nabla \cdot v_{n}, \, \xi_h)+
\sum_{T\in\Omega_h}\alpha_T (\nabla (p_n-\pi_h p_n),\nabla (\xi_h-\pi_h \xi_h))=0\quad\forall \xi_h\in L_h,
\]
and since $\nabla\cdot V_h\not\subset L_h$ only a perturbed divergence condition is satisfied
\begin{equation}
\label{eq:div:multigrid}
\|\nabla\cdot v_n\|_{L^2(\Omega)}=\|\nabla\cdot
(v-v_n)\|_{L^2(\Omega)}
\le c \|v-v_n\|_{H^1(\Omega)}.
\end{equation}
Fig.~\ref{fig:div} shows the total divergence $\|\nabla\cdot v(t)\|_{L^2(\Omega)}$ over the simulation domain for the coarse mesh solution (MG), the high fidelity solution (MG fine) and the DNN-MG solver for the experiments described in Sec.~\ref{sec:structure:numerics}. While convergence of the high fidelity is observed with respect to the coarse mesh solution, DNN-MG even further disturbs the divergence instead of improving it.

% \begin{figure}[t]
%   \captionsetup{format=plain}
%   \begin{captionbeside}{%
%       Divergence $\|\nabla\cdot v\|$ for the coarse mesh (MG),
%       the high fidelity (MG fine) and the deep neural network multigrid solution (DNN-MG).
%       Three results of the attempts to enforce divergence freedom are shown:
%       DNN-MG with the divergence of the DNN's output as penalty term (DNN-MG-$\mathcal{L}_{p_1}$),
%       DNN-MG with the divergence of the corrected velocity field as penalty term
%       (DNN-MG-$\mathcal{L}_{p_2}$) and enforcing divergence freedom by construction (DNN-MG-$\mathcal{L}_{\psi}$). \red{In my opinion it would be better to split this up. Have one figure with the old approaches here and then later this figure with the new ones. Space shouldn't be an issue}.\\ \bigskip
%     }[r]
%     \includegraphics[width=0.5\textwidth]{divergencefcts_convdiv.pdf}
%   \end{captionbeside}
%   \label{fig:div-old}
% \end{figure}

\subsection{Approaches for Enforcing Divergence Freedom}
\label{sec:structure:ansatzes}

In the following we will discuss three approaches for enforcing divergence freedom.
Numerical results will be presented afterwards in Sec.~\ref{sec:structure:numerics}.

\subsubsection{A Penalty Approach for Divergence Free Learning}
\label{sec:div12}

The first approach to enforce divergence freedom of the neural network output
of the DNN-MG method is through a modification of the loss function by a
squared divergence as penalty term. One then has 
\begin{equation}\label{loss:1}
  \mathcal{L}_{p_1}(v^L,v^{L+1}; d^{L+1}):= \mathcal{L}(v^L,v^{L+1};d^{L+1})
  +\gamma \sum_{n=1}^N
  \sum_{P_i\in\Omega_h} \big\Vert \nabla \cdot d^{L+1} \big\Vert^2_{L^2(P_i)},
\end{equation}
where $\mathcal{L}(v^L,v^{L+1};d^{L+1})$ is the original loss function in Eq.~\ref{loss} and $d^{L+1} = \mathcal{N}(v^L,\tilde{v}^{L+1})$ is the correction by the network based on the coarse solution $v_n^L$ and its prolongation $\tilde{v}^{L+1}$. The parameter $\gamma>0$ is used to control the weighting of the loss function.
Eq.~\ref{loss:1} penalizes the divergence of the correction and, through this, the neural network should at least not further distort the divergence of the velocity field.

An obvious extension of this idea is to penalize the divergence of the corrected solution $v^L+d^{L+1}$.
This is, in principle, advantageous since it could also correct for the divergence residual in the finite element formulation, cf.~Eq.~\ref{eq:div:multigrid}.
We thus introduce the second, modified loss function 
\begin{equation}\label{loss:2}
  \mathcal{L}_{p_2}(v^L,v^{L+1}; d^{L+1}):= l(v^L,v^{L+1};d^{L+1})
  +\gamma\sum_{n=1}^N \sum_{P_i\in\Omega_h} \big\Vert \nabla \cdot (v^L+d^{L+1}) \nabla\Vert^2_{L^2(P_i)} .
\end{equation}

\subsubsection{A Strictly Divergence Free Network Architecture}
\label{sec:div3}

Instead of including a penalty term in the loss function to ensure divergence freedom in a weak sense, our third approach modifies the neural network architecture so that only solenoidal vector fields are returned.
We accomplish this by learning the stream function $\psi$, which can be seen as a scalar potential for the divergence free velocity vector field~\cite{Chorin1993},
\begin{align}
  v = \nabla_{\perp} \, \psi = \begin{pmatrix}-\partial_y \\ \ \ \partial_x \end{pmatrix} \psi .
\end{align}
For an easy and efficient integration into the finite element framework, we are looking for divergence free velocity corrections $d^{L+1}$ by the neural network that are in the velocity function space, i.e.\ $d^{L+1}\in V_h^{L+1}$.
This holds when we learn stream functions $\psi_h \in L_h^{L+1}$ that are in the pressure space $L_h^{L+1}$ and the construction is then also again local on each patch.
%\[
%\psi_h = \begin{pmatrix}-\partial_y \\ \partial_x \end{pmatrix}\xi_h,
%\]
%where $\xi_h\in L_h$ are the functions of the pressure space.

For simplicity we assume that a patch $P_i$ corresponds exactly to one element on the fine mesh $\Omega_{L+1}$ and that it is a quadrilateral parallel to the axes on an integer lattice.
The 9 local scalar basis functions of the second order finite elements can then be defined in the simple monomial basis,
\[
L_h^{L+1}\big|_P = \big\{ \xi_{i=3\,\alpha+\beta} \equiv \xi_{\alpha,\beta}(x,y)=x^\alpha y^\beta,\; \alpha,\,\beta=0,\,1,\,2 \big\} ,
\]
and a stream function $\psi_h \in L_h^{L+1}$ is described by the $8$ basis function coefficients $s_{P_i}=(s_2,\dots,s_9)$ for this basis (the stream function has no harmonic part and $\xi_{0,0}$ can thus be disregarded).
%, compare the middle sketch in Fig.~\ref{fig:patch}.
Since we work with $L_h = W_h^{(2)}$, the velocity $v_h^{\psi}$ associated with $\psi_h$ lies in $L_h \vert_P$, as desired, and it can be represented in the eight velocity basis functions derived from the $\xi_{\alpha,\beta}(x,y)=x^\alpha y^\beta$,
\[
\begin{aligned}
  v_2^\psi&:=\nabla_\perp x = \begin{pmatrix}0\\1   \end{pmatrix},&\quad
  v_3^\psi&:=\nabla_\perp x^2= \begin{pmatrix}0\\2x   \end{pmatrix},&\quad
  v_4^\psi&:=\nabla_\perp y = \begin{pmatrix}-1\\0   \end{pmatrix}&\\
  v_5^\psi&:= \nabla_\perp xy = \begin{pmatrix}-x\\y   \end{pmatrix},&\quad
  v_6^\psi&:=\nabla_\perp x^2y = \begin{pmatrix}-x^2\\2xy   \end{pmatrix},&\quad
  v_7^\psi&:=\nabla_\perp y^2 = \begin{pmatrix}-2y\\0   \end{pmatrix}&\\
  v_8^\psi&:=\nabla_\perp xy^2 = \begin{pmatrix}-2xy\\y^2  \end{pmatrix},&\quad
  v_9^\psi&:=\nabla_\perp x^2y^2= \begin{pmatrix}-2x^2y\\2xy^2   \end{pmatrix}&\\
\end{aligned}
\]
By linearity we thus have
\begin{align}
  \label{eq:stream:velocity:basis}
  v_h^{\psi} \big\vert_{P_i} = \nabla_{\perp} \psi_h = \sum_{i=0}^8 s_i \, v_i^\psi
\end{align}
and the basis function coefficients for the stream function and the derived velocity are identical.

To be able to work within the usual finite element theory, we relate the non-standard velocity basis functions $v_i^\psi$ to pointwise values using the Lagrange basis
\begin{equation}
  \label{eq:Vh:Lagrange_basis}
V_h\big|_P=\{\phi_{1}^x,\dots,\,\phi_9^x,\;
\phi_{1}^y,\dots,\,\phi_9^y\}
\end{equation}
that satisfies
\[
\phi_i^x(x_j) =\begin{pmatrix}\delta_{ij}\\0\end{pmatrix},\quad
\phi_i^y(x_j) =\begin{pmatrix}0 \\\delta_{ij}\end{pmatrix},\quad i,j=1,\dots,9
\]
for the 9 Lagrange points $\delta_{ij}$.
Each of the derived velocity functions $v_i^\psi$  has a unique representation in this basis with coefficients $\eta_{ij}^x,\eta_{ij}^y$ for $i=2,\dots,9$ and $j=1,\dots,9$.
% such that
%\[
%v_i^{\psi} = \sum_{j=1}^9\eta_{ij}^x\phi_j^x + \eta_{ij}^y\phi_j^y
%\]
The pointwise values of the velocity vector field can thus be reconstructed as
\begin{equation}
  \label{eq:velocity:stream:pointwise}
  v_h^{\psi} \big\vert_{P_i} = \sum_{i=2}^9 s_i \sum_{j=1}^9 \eta_{ij}^x\phi_j^x + \eta_{ij}^y\phi_j^y .
\end{equation}
Along boundaries of the domain $\partial\Omega$ where Dirichlet
conditions are prescribed, the corresponding coefficients
$\eta_{ij}^x,\eta_{ij}^y$ are set to zero such that the resulting
update $v_h^\Psi$ is no longer strictly divergence free.

With the above, we can train for the eight stream function coefficients $s_{P_i}=(s_2,\dots,s_9)$ that through Eq.~\ref{eq:stream:velocity:basis} and Eq.~\ref{eq:velocity:stream:pointwise} correspond to a divergence free velocity correction $d^{L+1}(\{ s_i^{L+1} \} )$ satisfying $\|\nabla\cdot d^{L+1}\|_{L^2(\Omega)}=0$.
Similar to the DNN-MG setting, the corrections $d^{L+1}(\{ s_i^{L+1} \} )$ are averaged over all degrees of freedom shared by multiple patches.
For network training, we minimize
\[
\mathcal{L}_{\psi}(v^L,v^{L+1};s^{L+1}):=
\sum_{P_j\in \Omega_h} \big\Vert v^{L+1} - \big(v^L+ d^{L+1}(\{ s_i^{L+1} \} )\big) \big\Vert^2_{L^2(P_j)}.
\]

\subsection{Numerical Evaluation}
\label{sec:structure:numerics}

To evaluate the efficacy of the three approaches developed in Sec.~\ref{sec:structure:ansatzes} to improve the divergence freedom of the velocity vector field, we implemented these using PyTorch~\cite{PyTorch} and Gascoigne 3D~\cite{Gascoigne3d} based on the framework developed in~\cite{Margenberg2020}.
For the penalty-based approaches with the loss functions $\mathcal{L}_{p_1}$ and $\mathcal{L}_{p_2}$ the network and other aspects of the training process
remain unchanged compared to~\cite{Margenberg2020} where the network had 8634 parameters.
The network architecture we use for the stream function-based approach with $\mathcal{L}_{\psi}$ is also similar to the one in~\cite{Margenberg2020} but the convolutional layers are
replaced by fully connected ones of size $8\times 8$.
This results in a network with 9000 parameters.

\begin{figure}[t]
  \centering
  \begin{tikzpicture}[scale=5.5]
    \draw [fill=gray!20] (0,0) rectangle (2.25,0.4);
    \path[draw=black] (0,0) -- node[midway, above]{$\Gamma_{\textrm{wall}}$} ++(2.25,0) --
    node[midway, left]{$\Gamma_{\textrm{out}}$} ++(0,0.4) -- node[midway,
    below]{$\Gamma_{\textrm{wall}}$} ++(-2.25,0) -- node[midway, right]{$\Gamma_{\textrm{in}}$}cycle;
    \draw[fill=white] (0.25,0.1) rectangle (0.35, 0.2);
    \node at (1.55, 0.2){$\mathrm{Re}_{\text{test}}= 133$};
    \draw[decorate,decoration={brace,raise=1pt,amplitude=9pt,mirror}] (0,0) --
    node[below=9pt]{$2{.}25$} (2.25,0);
    \draw[decorate,decoration={brace,raise=1pt,amplitude=9pt}] (0,0) --
    node[above left=12pt and 12pt,rotate=90]{$0{.}4$} (0,0.4);
    \draw[decorate,decoration={brace,raise=1pt,amplitude=7pt,mirror}] (0.35,0.1) --
    node[right=7pt]{$0{.}1$} (0.35,0.2);
    \draw[decorate,decoration={brace,raise=1pt,amplitude=7pt}] (0.25,0.2) --
    node[above=7pt]{$0{.}1$} (0.35,0.2);
    % \draw[latex-latex] ([yshift=-2pt]0,0) -- node[fill=white]{$2{.}2$} ([yshift=-2pt]2.2,0);
    \draw[-latex] (0.55, 0.25) node[above]{$(0{.}3,\,0{.}15)$} --(0.3,0.15);
  \end{tikzpicture}
\caption{Geometry of the training scenario with a parabolic inflow profile $\Gamma_{\textrm{in}}$,
  do-nothing boundary conditions at the outflow boundary $\Gamma_{\textrm{out}}$ and no-slip conditions on the walls
  $\Gamma_{\textrm{wall}}$. The center of the obstacle is at $(0.3,\,0.15)$. For the test
  scenario the obstacle is shifted such that the center is at $(0.3,0.25)$.}
  \label{fig:geom}
\end{figure}
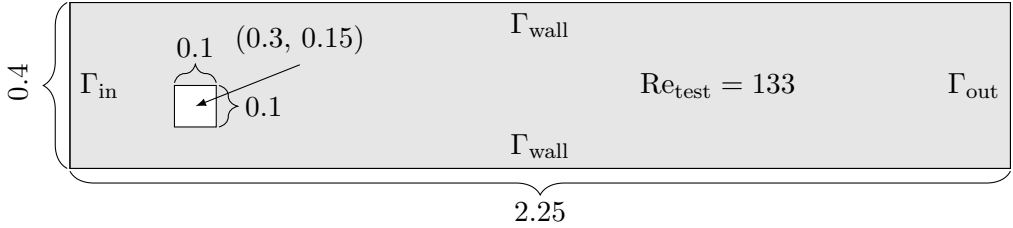

The flow we considered in our experiments is a modification of the classical benchmark of a flow around a cylinder~\cite{SchaeferTurek1996} where the obstacle is now a square, see Fig.~\ref{fig:geom}.
The modification simplified the implementation, since we were able to use uniform meshes, but does not change the qualitative behavior of interest to us and in particular one still has a non-stationary wake behind the obstacle as in the original benchmark~\cite{SchaeferTurek1996}.
The singularities at the corners of the obstacles have been excluded in the evaluation below.

\paragraph{Divergence freedom}
We first compare the divergence $\|\nabla\cdot v(t)\|_{L_2(\Omega_L)}$ obtained with the three approaches (DNN-MG-$\mathcal{L}_i$) as well as a coarse multigrid solution (MG), the original DNN-MG method proposed in~\cite{Margenberg2020} (DNN-MG), and, as reference, a fine multigrid solution on level $L+1$ (MG-fine).
The results are presented in Fig.~\ref{fig:div}.
%\begin{description}
%\item[MG] is the coarse mesh multigrid solution
%\item[MG (fine)] is the high fidelity  multigrid solution
%\item[DNN-MG] is the classical deep neural network multigrid approach without any special treatment of the divergence
%\item[DNN-MG-$\mathcal{L}_{p_1}$] is the penalty method based on penalyzing the network output (the correction) only, see~(\ref{loss:1})
%\item[DNN-MG-$\mathcal{L}_{p_2}$] is the penalty method based on penalyzing the corrected solution,see~(\ref{loss:2})
%\item[DNN-MG-$\mathcal{L}_{\psi}$] finally is based on the network architecture that guarantees strictly divergence free outputs, compare~Section~\ref{sec:div3}.
%\end{description}
%
They show that the approaches from Sec.~\ref{sec:structure:ansatzes} are not able to significantly reduce the divergence.
While DNN-MG-$\mathcal{L}_{\psi}$ based on the strictly divergence free output performs best, the reduction is still modest.
Ultimately, however, this is not surprising since
\[
\|\operatorname{div}\,d^{L+1}\|=0\quad\Rightarrow\quad
\|\operatorname{div}\,(v^L+d^{L+1})\|=
\|\operatorname{div}\,v^L\| ,
\]
and the perfectly divergence free correction does not address the divergence residual in $v^L$.
The second approach, using $\mathcal{L}_{p_2}$, was designed to address this but no reduction could be observed in practice.

We investigated the disappointing performance of the modified loss functions $\mathcal{L}_{p_1}$ and $\mathcal{L}_{p_2}$ that included a penalty term and
found that it stems from an unstable balance between the minimization of the divergence
and the fulfillment of the finite element Galerkin equation: If the weighting factor $\gamma$
for the divergence is chosen close to zero then the network behaves like the classical DNN-MG approach.
If we choose larger values for
$\gamma>0$, however, then the network is not able to sufficiently minimize the discrepancy to the high
fidelity Galerkin solution and overall the approximation quality suffers. We tested $\gamma$ in the
range $10^{-1},\cdots,\, 10^{-10}$ but did not find a value where the two error terms in
the loss functions $\mathcal{L}_{p_1}$ and $\mathcal{L}_{p_2}$ were sufficiently balanced.
For $\gamma > 10^{-7}$ we, in fact, observed that the penalty term causes the output to converge
towards zero, which is a divergence free velocity vector field but not the desired correction.
The effect was most pronounced for DNN-MG-$\mathcal{L}_{p_2}$.
%As expected the networks started to converge towards
%DNN-MG for even smaller $\gamma$.
%\textcolor{red}{Kann man etwas auf den Einfluss von $\gamma$ auf das Trainieren sagen?}

\begin{figure}[t]
  {\centering
    \includegraphics[width=0.5\textwidth]{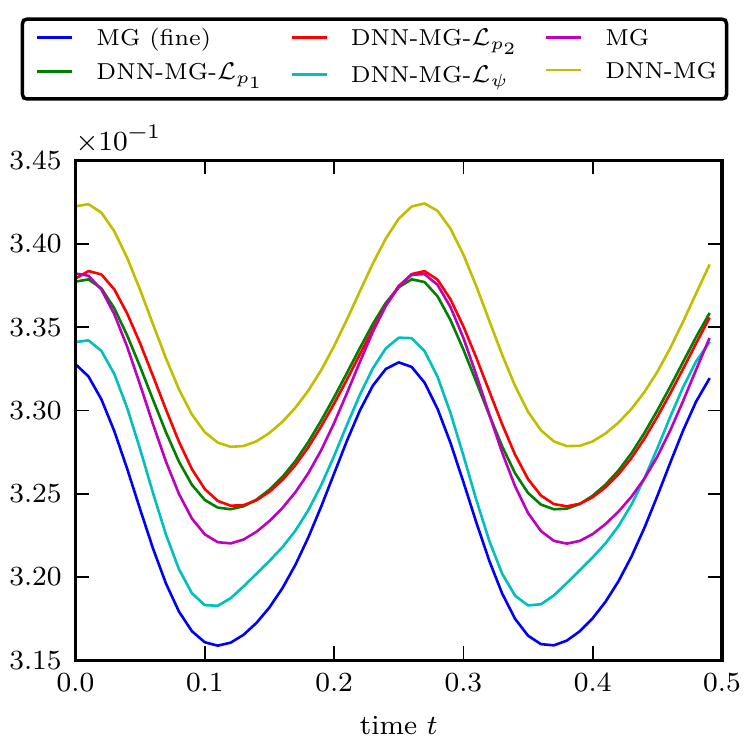}
    \includegraphics[width=0.5\textwidth]{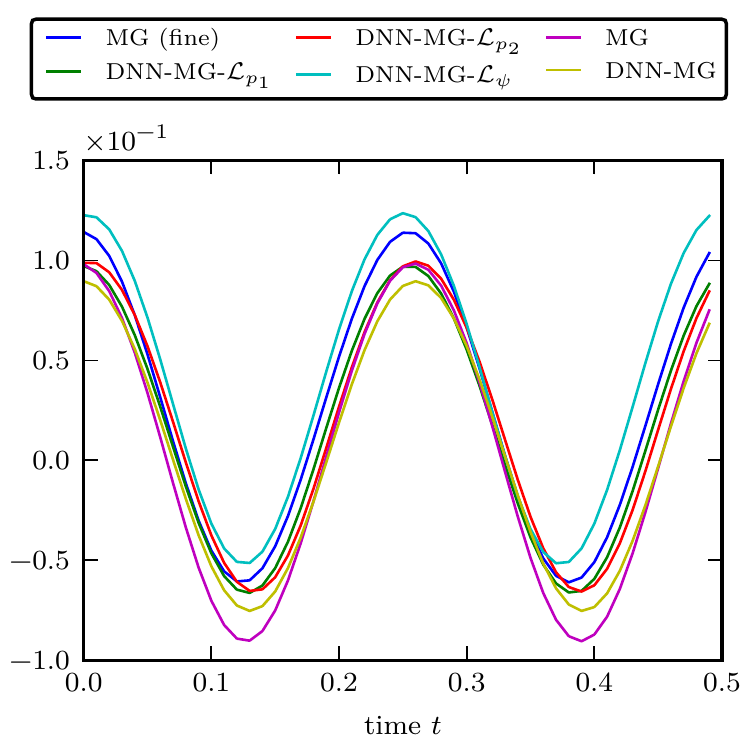}
  }
  \caption{Comparison of drag and lift functionals for the coarse mesh (MG),
    the high fidelity (MG fine) and the deep neural network multigrid solution (DNN-MG).
    Three results of the attempts to enforce divergence freedom are shown:
    DNN-MG with the divergence of the DNN's output as penalty term (DNN-MG-$\mathcal{L}_{p_1}$),
    DNN-MG with the divergence of the corrected velocity field as penalty term
    (DNN-MG-$\mathcal{L}_{p_2}$) and enforcing divergence freedom by construction (DNN-MG-$\mathcal{L}_{\psi}$).
  }
  \label{fig:draglift}
\end{figure}

\paragraph{Drag and Lift functionals}
Fig.~\ref{fig:draglift} shows the drag and lift functionals
of the obstacle for the different approaches. In
Table~\ref{tab:draglift} we also show the corresponding minimum and maximum
values, the mean values, and the amplitudes of the oscillation as well as
 the frequencies of the periodic flow pattern.
In the lower table of Table~\ref{tab:draglift} we report
the relative errors with respect to the high fidelity solution MG (fine).

The results show that the non-smooth corners of the
square obstacle have a negative effect on
the performance of the original DNN-MG approach since it no longer
provides the improvement in the functional values observed
in~\cite{Margenberg2020} for a smooth obstacle.
Among the approaches proposed in the present work,
DNN-MG-$\mathcal{L}_{\psi}$ based on learning the stream function
provides the most promising results.
It is able to reduce the relative
error in all functional outputs and almost matches the values of
MG (fine) in Fig.~\ref{fig:draglift}.
Furthermore, the training time for the neural network is in this
case only $\approx 75\%$ in comparison to the other DNN-MG approaches.

For DNN-MG-$\mathcal{L}_{\psi}$ we observe a deviation of the oscillation frequency
that is larger than for the other approaches.
We believe that this defect stems
from a complicated interplay of temporal dissipation and spatial
discretization. For a discussion we refer to \cite{Margenberg2020} and
also to~\cite{MargenbergRichter2020} where the same effect is
observed for a parallel time
stepping scheme.
In the context of the deep neural network multigrid
approach it must be investigated if this problem can be resolved by
training with high fidelity solutions that come from
discretizations that are finer in space and in time, i.e. by working with a finer time step.

\begin{table}[t]
  \begin{center} \footnotesize
    \begin{tabular}{l|cccc|cccc|c}
      \toprule
      &\multicolumn{4}{c|}{Drag}&\multicolumn{4}{c|}{Lift}\\
      Approach & Min & Max & Mean & Ampl. & Min & Max & Mean & Ampl. & Freq. \\
      \midrule
      coarse      & 0.3220 & 0.3382 & 0.3301 & 0.0162 & -0.0905 & 0.0985 & 0.0040 & 0.1890 & 3.8095 \\
      fine        & 0.3159 & 0.3329 & 0.3244 & 0.0170 & -0.0610 & 0.1144 & 0.0267 & 0.1754 & 3.9024 \\
      \midrule
      DNN-MG      & 0.3277 & 0.3425 & 0.3351 & 0.0148 & -0.0753 & 0.0896 & 0.0071 & 0.1650 & 3.8298 \\
      \midrule
      DNN-MG-$\mathcal{L}_{p_1}$ & 0.3240 & 0.3379 & 0.3309 & 0.0139 & -0.0665 & 0.0978 & 0.0156 & 0.1643 & 3.9130 \\
      DNN-MG-$\mathcal{L}_{p_2}$ & 0.3242 & 0.3384 & 0.3313 & 0.0142 & -0.0657 & 0.0999 & 0.0171 & 0.1657 & 3.9130 \\
      DNN-MG-$\mathcal{L}_{\psi}$ & 0.3178 & 0.3350 & 0.3264 & 0.0172 & -0.0544 & 0.1259 & 0.0357 & 0.1803 & 4.0179 \\
      \bottomrule
    \end{tabular}

    \bigskip
    \begin{tabular}{l|cccc|cccc|c}
      \toprule
      &\multicolumn{4}{c|}{Drag}&\multicolumn{4}{c|}{Lift}\\
      Approach & Min & Max & Mean & Ampl. & Min & Max & Mean & Ampl. & Freq. \\
      \midrule
      coarse      &1.93  &1.59  &1.76  &4.71  &48.4  &13.9  &85.0  &7.75  &2.38\\
      \midrule
      DNN-MG      &3.74  &2.88  &3.30  &12.9  &23.4  &21.7  &73.4  &5.93  &1.86\\
      \midrule
      DNN-MG-$\mathcal{L}_{p_1}$ &2.56  &1.50  &2.00  &18.2  &9.02  &14.5  &41.6  &6.33  &\bfseries 0.27\\
      DNN-MG-$\mathcal{L}_{p_2}$ &2.63  &1.65  &2.13  &16.5  &\bfseries 7.70  &12.7  &36.0  &5.53  &\bfseries 0.27\\
      DNN-MG-$\mathcal{L}_{\psi}$ &\bfseries 0.60  &\bfseries 0.63  &\bfseries 0.62  &\bfseries 1.18  &10.8  &\bfseries 10.1  &\bfseries 33.7  &\bfseries 2.79  &2.96\\
      \bottomrule
    \end{tabular}
  \end{center}
  \caption{Maximum, minumum values, mean and amplitude of oscillation
    for drag and lift functionals as well as frequency of the periodic
    solution. We indicate the results for the different finite element
    solution on coarse and fine mesh finite element simulation
    and for the different hybrid finite element / deep neural network
    approaches. While the upper table shows the raw results, the lower
    table indicates the relative error (in \%) with respect to the high fidelity
    solution. For each value we highlight the best result.}
  \label{tab:draglift}
\end{table}

\begin{figure}[t]
  \begin{center}
    \includegraphics[width=0.95\textwidth]{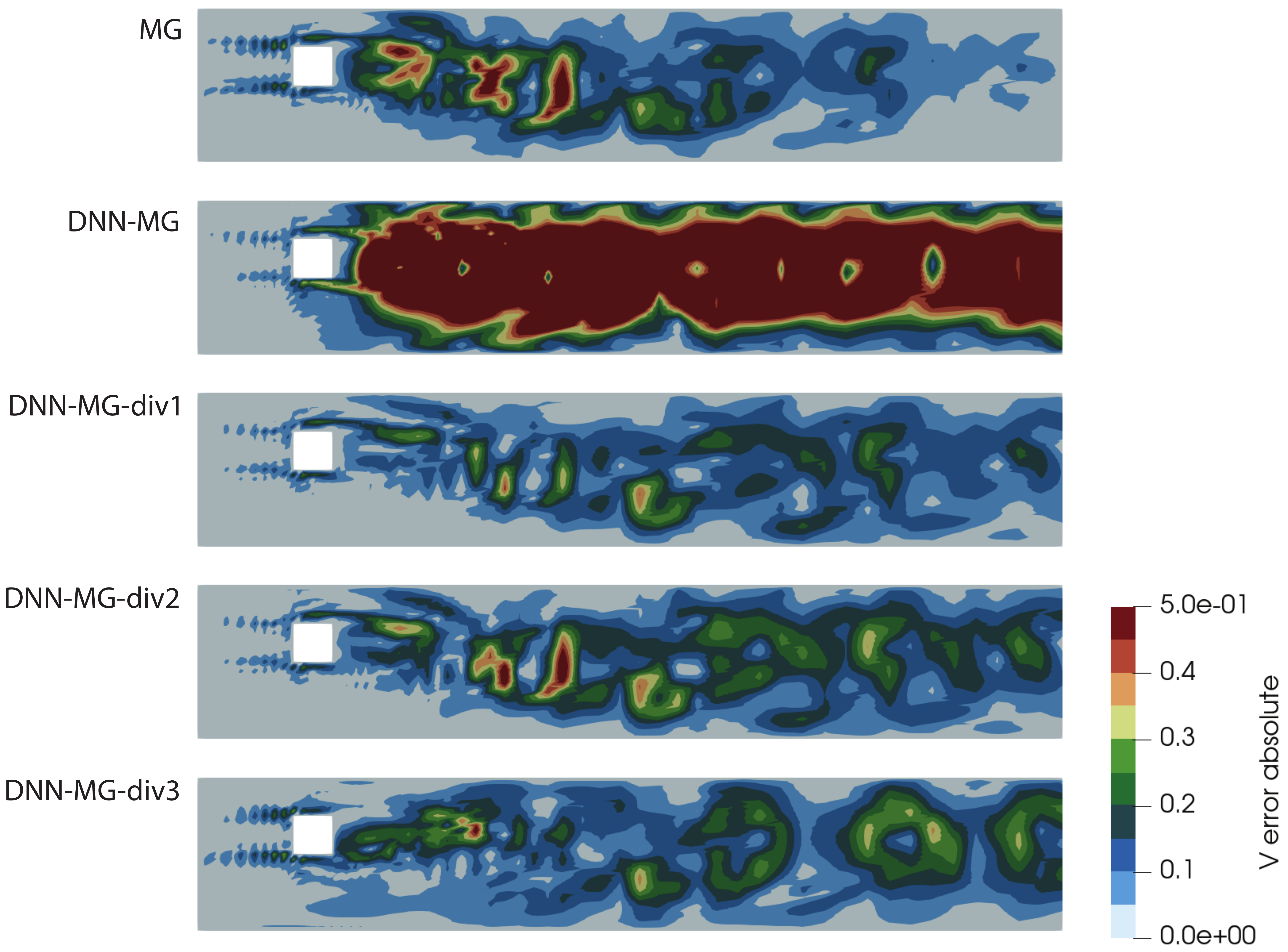}
  \end{center}
  \caption{Velocity error for the different approaches with respect to the high fidelity solution MG (fine).}
  \label{fig:velocity}
\end{figure}

\paragraph{Velocity vector field}
Fig.~\ref{fig:velocity} visualizes the error $|v(x)-v_{\mathrm{fine}}(x)|$ in the velocity vector field for MG, DNN-MG, DNN-MG-$\mathcal{L}_{p_1}$, DNN-MG-$\mathcal{L}_{p_2}$, DNN-MG-$\mathcal{L}_{\psi}$ with respect to the fine reference solution MG (fine).
To obtain meaningful comparisons, the outputs are synchronized such that we consider snapshots starting from a maximum of the drag functional, cf.~Fig.~\ref{fig:draglift}.
The results verify that all approaches proposed in Sec.~\ref{sec:structure:ansatzes} are able to reduce the error in the velocity field compared to MG and DNN-MG.

\paragraph{Discussion}
The results presented in this section are surprising in that the loss functions of Sec.~\ref{sec:structure:ansatzes} designed to reduce the divergence improve the overall quality of the obtained results, e.g.\ in terms of the lift and drag functionals, but do not lead to significantly better divergence freedom.
More work is hence required to obtain hybrid finite element / neural network-based simulations that improve this important property.
For the loss functions with the penalty factor it should in particular be reconsidered how these can be formulated so that their two constitutive terms do not compete with each other as we observed in our experiments.

The efficiency of DNN-MG-$\mathcal{L}_{\psi}$ based on the stream function is not entirely surprising to us.
The network there only has to learn a scalar function instead of two components of a vector field that, furthermore, satisfy a nontrivial constraint through the divergence freedom.
The stream function formulation also encodes the structure of the dynamics in a different way than the usual velocity one, cf.~\cite{Morrison1998a}, and this leads often to simulations that preserve structural properties such as energy or vorticity.
In the finite element context this was already observed by Fix~\cite{Fix1975}.
We believe that this partly also explains the reduced training time for DNN-MG-$\mathcal{L}_{\psi}$.

\section{Conclusion}
\label{sec:conclusion}

We have presented three approaches to enforce divergence freedom within the hybrid 
finite element / neural network DNN-MG method~\cite{Margenberg2020} for the simulation
of the Navier-Stokes equations. 
The first two of these use a penalty term that encourages divergence freedom while the third one
guarantees a divergence free correction by learning the stream function.
The penalty term-based approaches did not provide significant improvements compared to 
DNN-MG for both the divergence and the drag and lift functionals.
The third approach based on the stream function yielded significantly better overall 
results but also only a rather modest improvement for the divergence.

More work is hence required in the future to obtain hybrid 
finite element / neural network simulations that respect the divergence freedom of the Navier-Stokes equations.
For the penalty-based approaches it should be investigated how the two terms in the loss 
function can be balanced.
DNN-MG-$\mathcal{L}_2$ is thereby particularly interesting since it allows, in principle, to 
provide an overall reduction by correction $v^L$ towards a lower divergence.
For DNN-MG-$\mathcal{L}_{\psi}$ it would be interesting to investigate the efficiency when also the velocity-stream function formulation is used for the simulation. 

In future work, we also want to investigate the considerably more challenging problem of a hybrid finite element / neural network simulation for the Navier-Stokes equations in 3D. 
How other physical invariants, e.g.\ energy or Kelvin's circulation theorem, can be incorporated into hybrid simulations and what benefits this provides should also be investigated in more detail.
Stability, consistency and convergence are, furthermore, pertinent open questions that should be addressed.

\paragraph{Acknowledgement}
NM and TR acknowledge the financial support by the Federal Ministry of
Education and Research of Germany, grant number 05M16NMA as well as
the GRK 2297 MathCoRe, funded by the Deutsche Forschungsgemeinschaft,
grant number 314838170.
NM acknowledges support by the Helmholtz-Gesellschaft grant number HIDSS-0002 DASHH.
TR further acknowledges funding by the
Deutsche Forschungsgemeinschaft, grant number 411046898.
CL is funded by the Deutsche Forschungsgemeinschaft (DFG, German Research Foundation) – Project-ID 422037413 – TRR 287.

\bibliographystyle{plain}
\bibliography{../lit,neuralnets,../climate.bib,pubs.bib}

\end{document}